\documentclass[10pt,twocolumn,twoside]{ieeeconf}

\usepackage{cite}
\usepackage{amssymb,amsfonts,bm}
\usepackage{algorithmic}
\usepackage{graphicx}
\usepackage{textcomp}
\usepackage{amsmath,tikz}
\usepackage{color}
\usepackage{url}
\usepackage{mathtools}
\usepackage[hidelinks]{hyperref} 
\usepackage{rotating}
\usepackage{blkarray}
\usepackage{physics}
\usetikzlibrary{arrows}
\let\proof\relax
\let\endproof\relax
\usepackage{amsthm}
\newtheorem{thm}{Theorem}[section]

\newtheorem{rem}[thm]{Remark}

\usepackage{caption}
\usepackage{subcaption}
\usepackage{etoolbox}
\let\classAND\AND
\let\AND\relax
\usepackage{algorithm}
\usepackage{algorithmic}

\let\AND\classAND
\AtBeginEnvironment{algorithmic}{\let\AND\algoAND}

\floatname{algorithm}{Algorithm}

\graphicspath{{eps/}{png/}}
\DeclareSymbolFont{symbolsC}{U}{pxsyc}{m}{n}
\DeclareMathSymbol{\coloneqq}{\mathrel}{symbolsC}{"42}

\newcommand{\vertiii}[1]{{\left\vert\kern-0.25ex\left\vert\kern-0.25ex\left\vert
		#1 
		\right\vert\kern-0.25ex\right\vert\kern-0.25ex\right\vert}}
\usepackage{multirow}
\usepackage{makecell}
\usepackage{adjustbox}
\usepackage{float}

\DeclareMathOperator{\Diag}{diag}

\DeclareMathOperator{\const}{const}

\DeclareMathOperator{\vect}{vec}
\newcommand*{\vecbf}[1]{\mathbf{#1}} 
\newcommand*{\matbf}[1]{\mathbf{#1}} 
\newcommand*{\gbf}[1]{\bm{#1}} 
\newcommand*{\myset}[1]{\mathcal{#1}} 
\newcommand*{\lrangle}[1]{\langle #1 \rangle} 

\usepackage{pifont}
\usepackage{soul}

\begin{document}
	
	\title{DIRECT: A Differential Dynamic Programming Based Framework for Trajectory Generation}

	\author{Kun~Cao$^\ast$, Muqing~Cao$^\ast$, Shenghai~Yuan, and~Lihua~Xie,~\IEEEmembership{Fellow,~IEEE}
		\thanks{This research is supported by the National Research Foundation, Singapore under its Medium Sized Center for Advanced Robotics Technology Innovation. $^\ast$Equal contribution. K. Cao, M. Cao, S. Yuan and L. Xie (corresponding author) are with School 
			of Electrical and Electronic Engineering, Nanyang Technological University, 50 
			Nanyang Avenue, Singapore 639798 (email: kun001@e.ntu.edu.sg; \{mqcao, shyuan, elhxie\}@ntu.edu.sg).}}

	\maketitle
                           
\begin{abstract}
This paper introduces a differential dynamic programming (DDP) based framework for polynomial trajectory generation for differentially flat systems. In particular, instead of using a linear equation with increasing size to represent multiple polynomial segments as in literature, we take a new perspective from state-space representation such that the linear equation reduces to a finite horizon control system with a fixed state dimension and the required continuity conditions for consecutive polynomials are automatically satisfied. Consequently, the constrained trajectory generation problem (both with and without time optimization) can be converted to a discrete-time finite-horizon optimal control problem with inequality constraints, which can be approached by a recently developed interior-point DDP (IPDDP) algorithm. Furthermore, for unconstrained trajectory generation with preallocated time, we show that this problem is indeed a linear-quadratic tracking (LQT) problem (DDP algorithm with exact one iteration). All these algorithms enjoy linear complexity with respect to the number of segments. Both numerical comparisons with state-of-the-art methods and physical experiments are presented to verify and validate the effectiveness of our theoretical findings. The implementation code will be open-sourced\footnote{\url{https://github.com/ntu-caokun/DIRECT}}. 
\end{abstract}


\IEEEpeerreviewmaketitle

\section{Introduction}
\label{sec:intro}
Recent decade has witnessed many emerging applications of unmanned systems, such as search and rescue \cite{hu2013cooperative}, environmental monitoring \cite{ogren2004cooperative} and mapping \cite{yue2019multilevel}, etc. As one of the fundamental problems in robotics and control researches, trajectory planning of mobile robots has been actively studied, see \cite{mellinger2011minimum,bry2015aggressive, gao2020teach}.

Although polynomial-based approaches can generate energy-optimal smooth trajectories with fixed time allocation efficiently, they are usually formulated as a quadratic programming problem of big size and solved via various commercial solvers. 
In addition, although there are several methods considering bi-level optimization on energy and time in view of high nonlinearity and non-convexity of joint energy-time optimization, the computation complexity with respect to the number of segments is unknown in general. 
On the other hand, although DDP based methods enjoy linear complexity in the length of prediction horizon, they usually can only plan a relatively short time ahead for highly nonlinear systems such as quadrotor as the discretization step should be very fine.

\begin{figure}[t!]
    \centering
    \begin{subfigure}[b]{\linewidth}
        \includegraphics[width=0.98\linewidth]{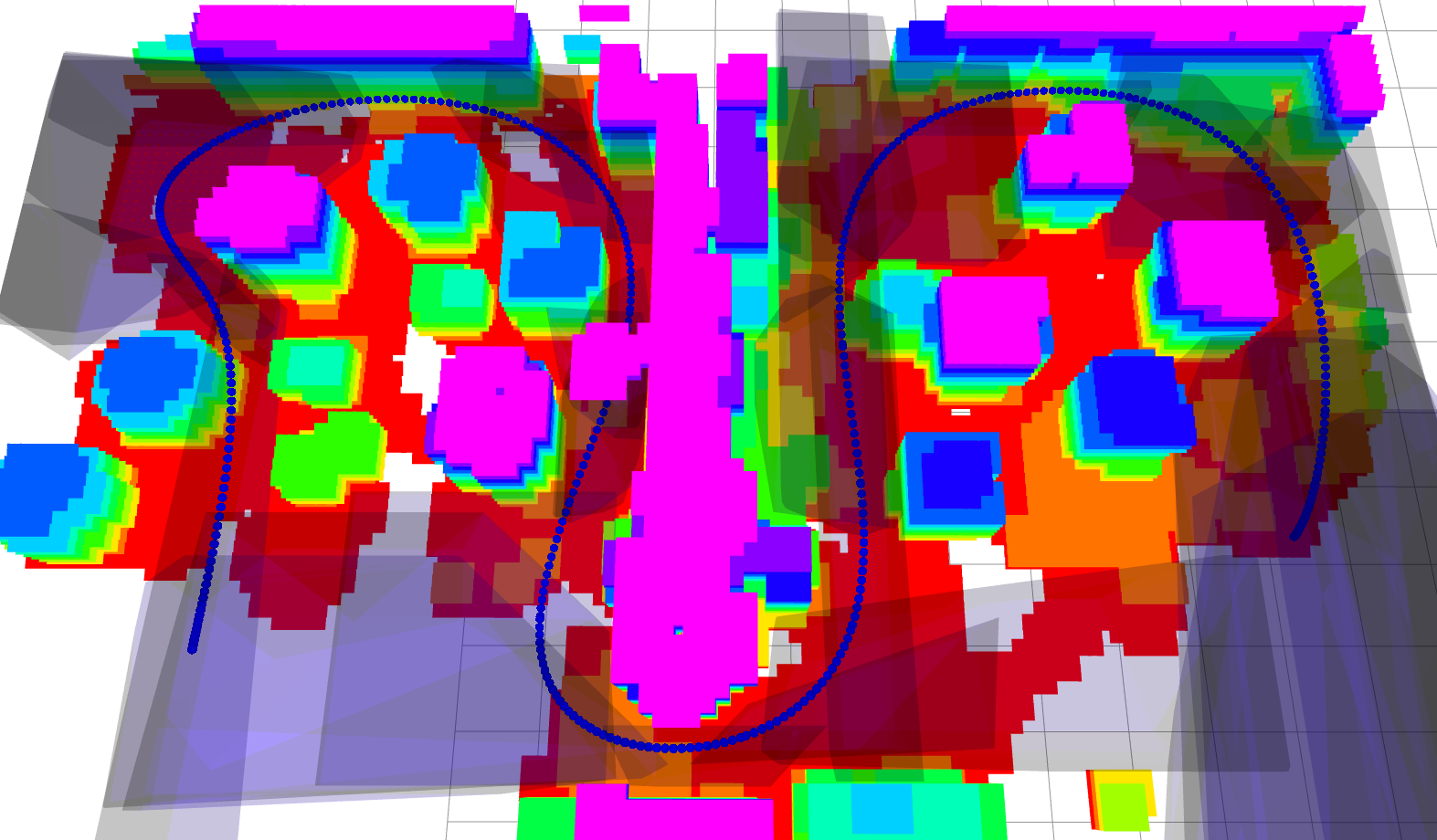}
    \end{subfigure}
    \par\bigskip    
    \centering
    \begin{subfigure}[b]{\linewidth}
        \includegraphics[width=0.98\linewidth]{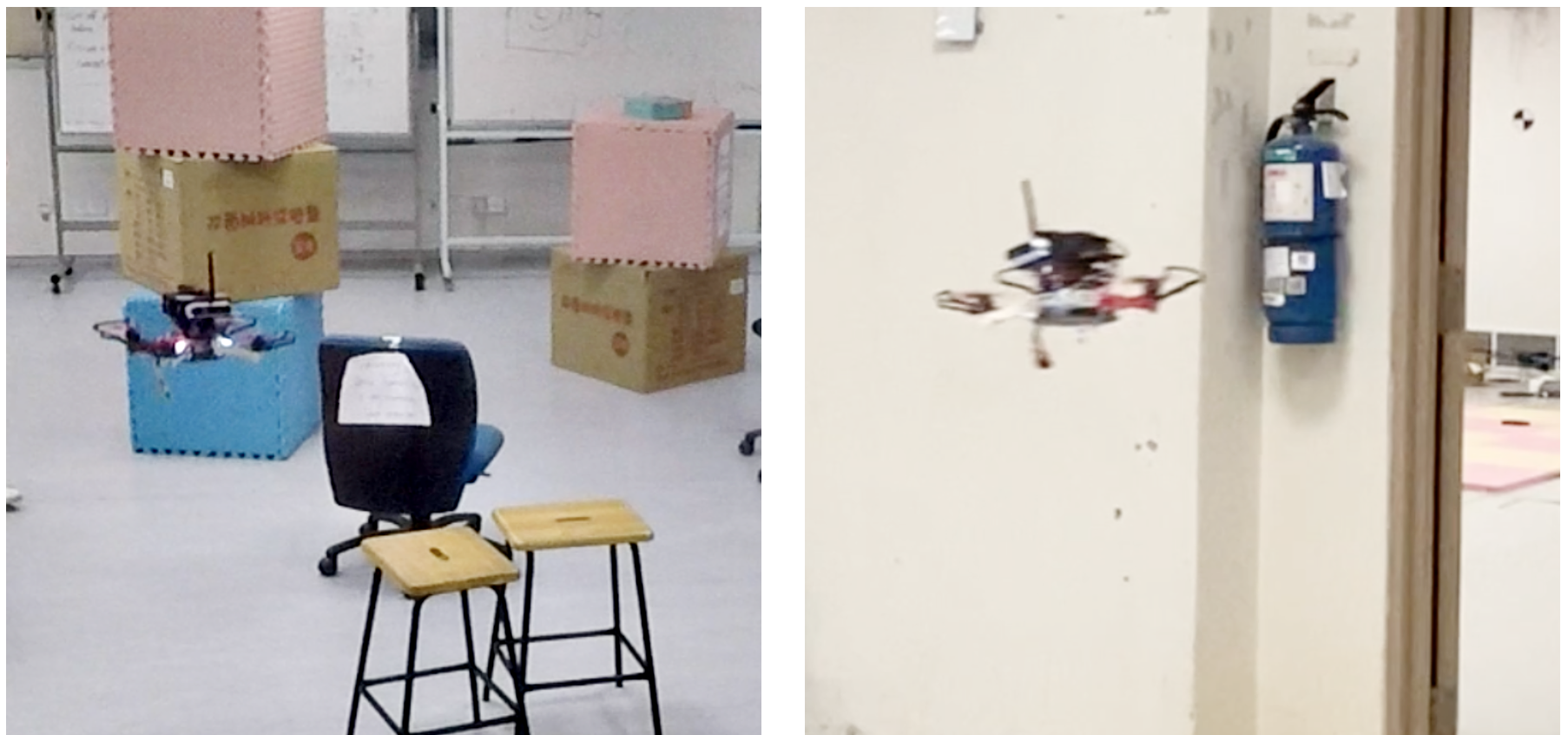}
    \end{subfigure}
    \caption{(Top) Trajectory generated using our proposed method in a two-room environment, the safe corridor consists of $34$ polyhedra. (Bottom) Snapshots of a quadrotor executing the trajectory generated.}
	\label{fig:longtraj}
\end{figure}

To address these issues, we propose DIRECT, a \underline{di}ffe\underline{re}ntial dynami\underline{c} programming based framework for \underline{t}rajectory generation. In particular, this framework uses a state-space equation to characterize the polynomial trajectory such that half of its coefficients can be saved and the continuity conditions on the derivatives of two consecutive segments can be automatically satisfied. By leveraging the differential flatness of quadrotor systems and this perspective of polynomial trajectory, we further reformulate the trajectory generation problem into a finite-horizon free-end discrete-time optimal control problem such that it can be addressed by DDP algorithm without fine discretization and the help of solvers. Unlike many existing approaches where optimization with respect to time and polynomial coefficients are carried out alternately, this algorithm features joint optimization of all variables directly, which explains the naming of this framework.

The proposed approach is based on the observation that the state at the proximal endpoint of an odd order polynomial is solely determined by half of the coefficients (low-order terms) while the state at the distant endpoint is co-determined by all of the coefficients and time. Along with the continuity condition, a long trajectory of multiple segments can be characterized by the initial state and a sequence of control inputs which are functions of the other half of the coefficients (high-order terms) and time. 
By incorporating the final state constraint into the terminal cost and introducing the safety and dynamical feasibility constraints on the state and control variables via control points, the optimization problem is reformulated into a discrete-time finite-horizon optimal control problem with inequality constraints. The resulting problem either with or without time optimization can be solved under the same framework via a recently developed IPDDP algorithm in \cite{pavlov2021interior} with linear complexity. In addition, we show that this problem can be further simplified to a LQT problem if there is no constraint and hence admits an analytical solution, which is still within the DDP framework (exactly one round of backward forward iteration). Our contributions can be summarized as follows:
\begin{itemize}
	\item We show that an odd-order piecewise-polynomial trajectory can be characterized by state-space equations;
	\item We propose a DDP-based framework to generate polynomial trajectories with / without time optimization and constraints;
	\item We achieve comparable results with state-of-the-art methods and perform extensive real-world experiments. The code will be open-sourced.
\end{itemize}
 
The rest of this paper is organized as follows. Some related works are presented in Section \ref{sec:related}. Section \ref{sec:ddp} details the DDP-based framework for trajectory generation. Some specializations of this scheme are discussed in Section \ref{sec:spec}. Section \ref{sec:results} gives some comparisons with state-of-the-art methods and real-world implementations. Section \ref{sec:conclu} draws the conclusion.

\emph{Notations:} In this paper, the sets of real, positive real and nonnegative real numbers are denoted by $\mathbb{R}$, $\mathbb{R}_{+}$ and $\mathbb{R}_{\geq 0}$, respectively. Denote by $\mathbb{R}^{n}$ and $\mathbb{R}^{m \times n}$ the sets of $n$-dimensional real vectors and $m \times n $ real matrices, respectively. Let $\mathbb{Z}_{+}$ and $\mathbb{Z}_{\geq 0}$ be the sets of positive and nonnegative integers, respectively. Let $\|\vecbf{x}\|$  be the $2$ norm of $\vecbf{x}$ and $\|\vecbf{x}\|_{\matbf{A}}^{2} = \vecbf{x}^{\top}\matbf{A}\vecbf{x}$. Denote by $\matbf{A}^{\top}$ and $\matbf{A}^{-1}$ the transpose and inverse of $\matbf{A} \in \mathbb{R} ^{n \times n}$, respectively. Let $\matbf{I}_{n} \in \mathbb{R}^{n \times n}$ be the $n$-dimensional identity matrix and $\vecbf{1}_{n}$ the $n$-dimensional column vector with all entries of $1$. Let $\otimes$ denote the Kronecker product and $\myset{I}_{n} = \{0,\dots,n-1\}$. Let $\vecbf{x}^{(i)}$ be the $i$-th order derivative of vector $\vecbf{x}$ and $\vecbf{x}^{[i]} = [\vecbf{x}^{(0)}, \dots, \vecbf{x}^{(i-1)}]$. Denote the vectorization operation by $\vect(\cdot)$, i.e., $\vect([\vecbf{a}, \vecbf{b}]) = [\vecbf{a}^{\top}, \vecbf{b}^{\top}]^{\top}$.

\section{Related Works}
\label{sec:related}
\subsection{Trajectory planning for quadrotors}
The seminal work \cite{mellinger2011minimum} firstly proved the differential flatness for quadrotor and proposed the minimum snap piecewise polynomial trajectory generation problem, which can be transformed into a quadratic programming (QP) problem with safety corridors being included as linear constraints on discretized points and hence can be solved efficiently via solver. The time allocation and dynamical feasibility problems were approached by projected gradient with backtracking line search in an outer loop and temporal scaling, respectively. Later, the authors in \cite{bry2015aggressive} proposed a closed-form solution for the trajectory generation problem with a generalized cost where the time allocation, safety, and dynamical feasibility aspects are handled by gradient descent of joint energy-time cost, addition of waypoints, and temporal scaling, respectively. In \cite{liu2017planning}, the authors presented an efficient safe flight corridor construction method and a temporal scaling step (similar to \cite{bry2015aggressive}) on time allocation to generate safe and dynamical feasible trajectories. To avoid high nonlinearity related to joint space-time optimization, some recent works propose methods that iteratively optimize spatial and time variables in two separate subroutines. The authors in \cite{gao2020teach} introduced a spatial-temporal optimization method which iteratively refines the geometrical coefficients and time-warping functions via QP and second order cone program (SOCP), respectively. In \cite{sun2021fast}, the authors decomposed the trajectory optimization problem as a bi-level optimization problem where the gradient required in the outer loop is analytically obtained from the dual solution of the inner loop QP. An alternating minimization scheme was proposed in \cite{wang2020alternating} where complex safety constraints were not considered. Recently in \cite{WANG2021GCOPTER}, the authors formulate an unconstrained nonlinear optimization that can be solved by quasi-Newton methods efficiently, but the safety and dynamical feasibility are enforced as penalty functions on discretized points instead of hard constraints.
In terms of large scale trajectory ($\geq 10^{3}$), the authors in \cite{burke2020generating} proved that the minimum snap problem with fixed time allocation can be solved in linear complexity by exploiting the structure of optimality conditions. The computational complexity has been further improved in \cite{wang2020generating} by eliminating a matrix inversion operation. An outer loop with analytical projected gradient was introduced for \cite{burke2020generating} in \cite{burke2021fast} such that the time allocation can be optimized. However, the computational complexity is unspecified for the case with safety and dynamical feasibility constraints. 

\subsection{DDP algorithms}
The DDP algorithm firstly introduced in \cite{mayne1966second} enjoys linear computational complexity (w.r.t. horizon) and local quadratic convergence \cite{de1988differential} and has been applied in trajectory optimization of  various systems, e.g., biped walking \cite{morimoto2003minimax} and robotic manipulation \cite{kumar2016optimal}. Much research has been devoted to generalizing the DDP to the case with inequality constraints and can be generally classified into two categories. The first category converts the constrained problems to unconstrained ones via penalty \cite{plancher2017constrained} while the other deals with the constraints explicitly by identifying the active ones \cite{xie2017differential}. To avoid the combinatorial problem regarding the active constraints, a constrained version of Bellman's principle of optimality has been introduced in \cite{pavlov2021interior,aoyama2020constrained}, which augments the control input with dual variables. However, these DDP algorithms have only been applied to a short-duration (e.g. $2\sim3\mathrm{s}$, $200\sim300$ steps of sampling period $0.01 \mathrm{s}$) quadrotor planning problem and takes quite long time for simple constraints ($\geq 5 \mathrm{s}$, see Table 7 in \cite{aoyama2020constrained}). In this paper, we shall exploit the differential flatness of quadrotors and the linear complexity property of DDP algorithm to generate long-duration smooth trajectories efficiently.

\section{DDP based framework for trajectory generation}
\label{sec:ddp}
In this section, we shall first introduce our proposed state-space representation for polynomial trajectories and then present our algorithm for jointly optimizing energy and time. 

\subsection{State-space representation}
We consider the problem of generating an $N$-segment $n$-th order polynomial trajectory in $d$-dimensional space with the $k$-th segment $\vecbf{p}_{k} \in \mathbb{R}^{d}$ being parametrized by:
\begin{equation}
\label{eq:single_seg}
	\vecbf{p}_{k}(t) = \matbf{C}_{k}^{\top}\vecbf{b}(t), t \in [0, t_{k}],
\end{equation}
where $\matbf{C}_{k} \in \mathbb{R}^{(n+1)\times d}$ is the coefficient matrix and $\vecbf{b}(t) = [1, t, \dots, t^{n}]^{\top}$ is the $n$-order monomial basis. Hence, the entire trajectory on $[0, \tau_{N-1}]$ can be defined with $$\vecbf{p}(t) = \vecbf{p}_{k}(t - \tau_{k-1}) = \matbf{C}_{k}^{\top}\vecbf{b}(t - \tau_{k-1})$$ for $t \in [\tau_{k-1}, \tau_{k}]$ where $\tau_{-1}=0$ and $\tau_{k} = \sum_{i=0}^{k} t_{k}$, $k \in \myset{I}_{N}$. 

Following \cite{mellinger2011minimum,wang2020generating,WANG2021GCOPTER}, we consider the problem of enforcing the smoothness of trajectory up to the $(m-1)$-th ($m = \frac{n+1}{2}$) order derivative and optimizing on the energy cost up to the $m$-th order derivative. 

We consider a single segment $k$. By \eqref{eq:single_seg}, the collection of trajectory up to $(m-1)$-th order derivative can be written as follows:
\begin{equation}
\label{eq:pk_derivative}
	\vecbf{p}_{k}^{[m]}(t) = \matbf{C}_{k}^{\top}\vecbf{b}^{[m]}(t).
\end{equation}
Collecting both the endpoints, one has:
\begin{equation}
 	[\vecbf{p}_{k}^{[m]}(0), \vecbf{p}_{k}^{[m]}(t_{k})]
		  =  \matbf{C}_{k}^{\top}[\vecbf{b}^{[m]}(0), \vecbf{b}^{[m]}(t_{k})] := \matbf{C}_{k}^{\top}\matbf{F}_{k}^{\top},
\end{equation}
where $\matbf{F}_{k}$ can be expressed explicitly with a $2$-by-$2$ block matrix \cite{wang2020generating}, i.e., $\matbf{F}_{k} = \left[\begin{smallmatrix}
\matbf{F}_{k,11} & \mathbf{0} \\
\matbf{F}_{k,21} & \matbf{F}_{k,22} \\
\end{smallmatrix} \right]$, where each block is of size $m \times m$, and 
\begin{equation}
	\begin{aligned}
		[\matbf{F}_{k,11}]_{ij} &= \begin{cases}
		\begin{aligned}
			(i-1)!, &~\mathrm{if}~ i=j, \\
			0, &~\mathrm{otherwise},
		\end{aligned}
		\end{cases}	\\
		[\matbf{F}_{k,21}]_{ij} &= \begin{cases}
		\begin{aligned}
		\frac{(j-1)}{(j-i)!}t_{k}^{j-i}, &~\mathrm{if}~ i\leq j, \\
		0, &~\mathrm{otherwise},
		\end{aligned}
		\end{cases} \\
		 [\matbf{F}_{k,22}]_{ij} &=
		\frac{(m+j-1)!}{(m+j-i)!}t_{k}^{m+j-i}.	
	\end{aligned}
\end{equation}
Next, we partition $\matbf{C}_{k}$ according to the block structure of $\matbf{F}_{k}$, i.e., $\matbf{C}_{k} = [\matbf{C}_{k,1}^{\top}, \matbf{C}_{k,2}^{\top}]^{\top}$ and obtain the following equality:
\begin{equation}
\label{eq:ss_eq}
	(\vecbf{p}_{k}^{[m]}(t_{k}))^{\top} = \matbf{F}_{k,21}\matbf{F}_{k,11}^{-1} (\vecbf{p}_{k}^{[m]}(0))^{\top} + \matbf{F}_{k,22} \matbf{C}_{k,2},
\end{equation} 
where $\matbf{C}_{k,1}$ was substituted with $\matbf{F}_{k,11}^{-1} (\vecbf{p}_{k}^{[m]}(0))^{\top}$. 

\begin{figure}
	\centering
	\includegraphics[scale = 0.6,trim={1cm 0cm 1cm 0cm},clip]{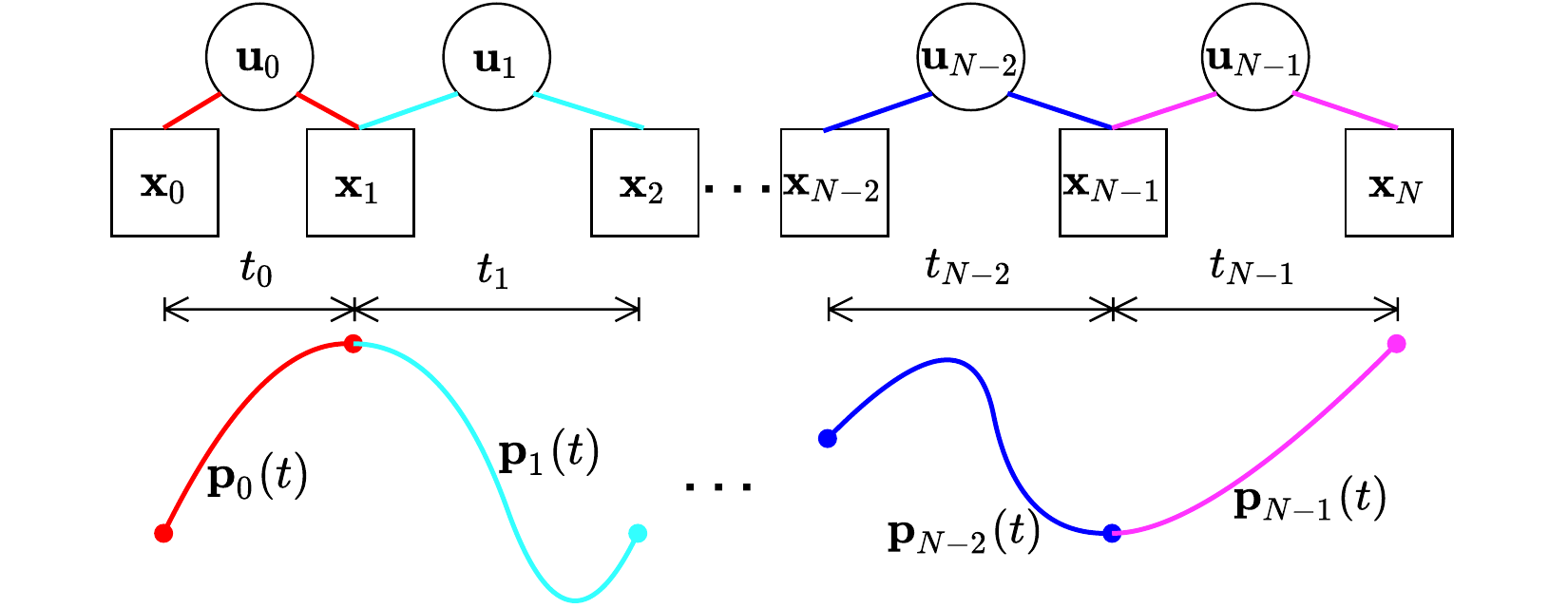}
	\caption{State space representation.}
	\label{fig:state_space}
\end{figure}

Letting $\vecbf{x}_{k} = \vect(\vecbf{p}_{k}^{[m]}(0))$, $\vecbf{v}_{k} = \vect(\matbf{C}_{k,2}^{\top})$, $\vecbf{u}_{k} = [\vecbf{v}_{k}^{\top}, t_{k}]^{\top}$ and by the continuity condition $\vecbf{p}_{k}^{[m]}(t_{k}) = \vecbf{p}_{k+1}^{[m]}(0)$, one has 
\begin{equation}
\label{eq:state}
\begin{aligned}
\vecbf{x}_{k+1} = \matbf{A}(t_{k}) \vecbf{x}_{k} + \matbf{B}(t_{k})\vecbf{v}_{k} := \vecbf{f}(\vecbf{x}_{k}, \vecbf{u}_{k}),
\end{aligned}
\end{equation}
where $\matbf{A}(t_{k}) =   (\matbf{F}_{k,21}\matbf{F}_{k,11}^{-1}) \otimes \matbf{I}_{d}$ and  $\matbf{B}(t_{k}) = \matbf{F}_{k,22} \otimes \matbf{I}_{d}$. As shown in Fig. \ref{fig:state_space}, this operation transfers the original representation of trajectory into a state-space equation, where the system state is the state (position, velocity, acceleration, etc.) at each endpoint and the control input is a concatenation of half of the coefficients of the segment (i.e., $\matbf{C}_{k,2}$) and duration (i.e., $t_{k}$). Unlike previous application of DDP algorithms in quadrotor planning and control where direct discretization is used, this representation enables the characterization of a long-duration trajectory with fewer number of variables. 

\subsection{Objective function}
Previous results on generating polynomial trajectories (e.g., \cite{mellinger2011minimum,wang2020generating}) usually consider the so-called interpolating splines \cite{zhang1997splines}, where the spline must pass some prescribed waypoints exactly at specific times. However, for the trajectory generation problem for time-critical aggressive flights, quite often it is sufficient that the generated trajectories only pass through the vicinity of these waypoints, i.e., smoothing splines \cite{sun2000control}. In this paper, we aim to find a smoothing spline such that the following objective function can be minimized:
\begin{equation}
\label{eq:obj0}
\begin{aligned}
	 J(\myset{C}, \myset{T}) =  &\sum_{k=0}^{N-1}(  \|\vecbf{x}_{k} - \vecbf{x}_{k,\mathrm{g}}\|_{\matbf{Q}_{k}}^{2} + \int_{0}^{t_{k}} \sum_{i=1}^{m} \eta_{i,k} \|\vecbf{p}_{k}^{(i)}(t)\|^{2}\dd t \\ &\quad \quad + w_{k} t_{k}^{2})  + \|\vecbf{x}_{N} - \vecbf{x}_{N,\mathrm{g}}\|_{\matbf{Q}_{N}}^{2},
\end{aligned}
\end{equation}
where $\myset{C} := \{\matbf{C}_{k}\}_{k \in \myset{I}_{N}}$ and $\myset{T} := \{t_{k}\}_{k \in \myset{I}_{N}}$.
In this cost function, the first term is designed to minimize the discrepancy between the current state $\vecbf{x}_{k}$ and the desired state $\vecbf{x}_{k,\mathrm{g}}$ at the beginning of each segment, which ``attracts" the endpoint to a prescribed waypoint; the second and third terms penalize the aggressiveness of planned trajectory up to the $m$-th order and flight time during the $k$-th segment, respectively; the last term penalizes the discrepancy between the terminal state $\vecbf{x}_{N}$ and the goal state $\vecbf{x}_{N,\mathrm{g}}$; and $\matbf{Q}_{k}$, $\eta_{i,k}$, $w_{k}$, ($k \in \myset{I}_{N}$), and $\matbf{Q}_{N}$ denote the weights for each term. It can be found that different from the literature (e.g., \cite{mellinger2011minimum,wang2020generating}) where hard constraints for the endpoints of $\vecbf{p}_{k}$ were used, we have softened these constraints by incorporating them into an objective function. 

For usually considered minimum jerk or minimum snap trajectories, $\matbf{Q}_{k} = \matbf{0}$, $\eta_{i,k} = 0$ for all $i = 1, \dots, m-1$ and $k \in \myset{I}_{N}$. In this case, one can partition $\vecbf{b}^{(m)}(t) (\vecbf{b}^{(m)}(t))^{\top}$ into $2$-by-$2$ block matrix with the size of each block being $m \times m$. 

Observing that only the right-bottom block is non-zero and denoting it by $\matbf{R}_{\vecbf{v}}(t)$, one can rewrite \eqref{eq:obj0} as
\begin{equation}
\label{eq:obj1}
	J(\myset{U}) =\sum_{k=0}^{N-1}(  \|\vecbf{x}_{k} - \vecbf{x}_{k,\mathrm{g}}\|_{\matbf{Q}_{k}}^{2} + \|\vecbf{u}_{k}\|_{\matbf{R}_{k}}^{2})   + \|\vecbf{x}_{N} - \vecbf{x}_{N,\mathrm{g}}\|_{\matbf{Q}_{N}}^{2}
\end{equation} with $\myset{U} := \{\vecbf{u}_{k}\}_{k \in \myset{I}_{N}}$ and
 $\matbf{R}_{k} = \Diag(\eta_{m, k} \int_{0}^{t_{k}} \matbf{R}_{\vecbf{v}}(t) \dd t, w_{k})$.

\begin{rem}
For the case with $\matbf{Q}_{k} \neq \matbf{0}$, $\eta_{i,k} > 0$ for some $i = 1, \dots, m-1$ and $k \in \myset{I}_{N}$, an equation similar to \eqref{eq:obj1}  can be obtained by using \begin{equation*}
	\begin{aligned}
	\|\vecbf{p}_{k}^{(i)}(t)\|^{2} & = (\vecbf{b}^{(i)}(t))^{\top} \matbf{C}_{k} \matbf{C}_{k}^{\top}\vecbf{b}^{(i)}(t) \\
	& = (\vect(\matbf{C}_{k}^{\top}))^{\top} [\vecbf{b}^{(i)}(t) (\vecbf{b}^{(i)}(t))^{\top} \otimes \matbf{I}_{d}] \vect(\matbf{C}_{k}^{\top})
	\end{aligned}
\end{equation*} and 
$\vect(\matbf{C}_{k}^{\top}) =  \left[\begin{smallmatrix}
	(\matbf{F}_{k,11}^{-1}\otimes \matbf{I}_{d}) \vecbf{x}_{k} \\ \vecbf{v}_{k}
\end{smallmatrix} \right]$. In addition, while other general differentiable functions can be used for each cost term in \eqref{eq:obj0}, we use a quadratic form for simplicity.
\end{rem}

\subsection{Constraints}
In general, there are two types of methods handling the safety and dynamical feasibility constraints. The first method discretizes each segment of trajectory and adds inequality constraints on the states at these time instants, while the issue is that the feasibility for the continuous time interval cannot be guaranteed. The other method is to add constraints on a set of control points, which are the vertices of the convex hull containing the segment, but it suffers from conservativeness. A recent work \cite{tordesillas2020minvo} proposes the MINVO basis, which has been shown to be less conservative than the widely used Bezier basis.  We use the MINVO basis to construct a set of control points for our polynomial trajectory. Denote the $i$-th order control points by $\gbf{\Xi}_{k}^{\langle i \rangle} \in \mathbb{R}^{(n+1-i) \times d}$, which can be computed by:
\begin{equation}
\label{eq:control_points}
	\gbf{\Xi}_{k}^{\lrangle{i}} = \matbf{T}^{\lrangle{i}}(t_{k}) \matbf{C}_{k},
\end{equation}
where $\matbf{T}^{\lrangle{i}}(t_{k}) \in \mathbb{R}^{(n+1-i) \times (n+1)}$ represents the transformation matrix from polynomial coefficients to the $i$-th order control points. 

Assuming that a safe flight corridor (consists of $N$ convex polyhedra) can be generated from perception system and polyhedron $k$ is represented by $s_{k}$ hyperplanes. With \eqref{eq:control_points}, the safety constraint can be written as:
\begin{equation}
\label{eq:cons_1}
	\matbf{W}_{k} \cdot (\gbf{\Xi}_{k}^{\lrangle{0}})^{\top} \preceq (\vecbf{1}_{n+1}^{\top} \otimes \vecbf{h}_{k}^{\lrangle{0}}),
\end{equation}
where each row of $\matbf{W}_{k} \in \mathbb{R}^{s_{k} \times d}$ and $\vecbf{h}_{k}^{\lrangle{0}} \in \mathbb{R}^{s_{k}}$ corresponds to the parameters of a single hyperplane in polyhedron $k$. 

On the other hand, the dynamical feasibility is encoded by the following box constraint on the $i$-th order control points $\gbf{\Xi}_{k}^{\langle i \rangle}$ with $i \in [1,m-1]$:
\begin{equation}
\label{eq:cons_2}
 \underline{h}_{k}^{\lrangle{i}} \matbf{I}_{d \times (n+1-i)} \preceq	 (\gbf{\Xi}_{k}^{\lrangle{i}})^{\top} \preceq \overline{h}_{k}^{\lrangle{i}} \matbf{I}_{d \times (n+1-i)},
\end{equation} where $\underline{h}_{k}^{\lrangle{i}}  < 0$ and $\overline{h}_{k}^{\lrangle{i}} > 0$ are the lower and upper bounds (e.g., minimal / maximal velocity / acceleration), respectively.

Additionally, one more inequality constraint on time can be introduced to avoid negative time allocation during optimization:
\begin{equation}
\label{eq:time_ineq}
    t_{k} \geq \underline{t},
\end{equation}
 where $\underline{t} \in \mathbb{R}_{+}$. 

Concatenating the vectorized forms of \eqref{eq:cons_1} and \eqref{eq:cons_2} with \eqref{eq:time_ineq}, one has
\begin{equation}
\label{eq:cons}
	\vecbf{g}_{k} := \Diag(\matbf{H}_{k}, -1) [
	 \vecbf{x}_{k}^{\top}, \vecbf{u}_{k}^{\top}]^{\top} - [\vecbf{h}_{k}^{\top}, -\underline{t}]^{\top} \preceq \vecbf{0},
\end{equation} where 
\begin{equation*}
	\matbf{H}_{k} = \left[\begin{smallmatrix}
	 \matbf{T}^{\lrangle{0}}(t_{k}) \otimes \matbf{W}_{k} \\
	 - \matbf{T}^{\lrangle{1}}(t_{k}) \otimes \matbf{I}_{d} \\
	 \matbf{T}^{\lrangle{1}}(t_{k}) \otimes \matbf{I}_{d} \\
	 \vdots \\
	 -\matbf{T}^{\lrangle{m-1}}(t_{k}) \otimes \matbf{I}_{d} \\
\matbf{T}^{\lrangle{m-1}}(t_{k}) \otimes \matbf{I}_{d} \\	 
	\end{smallmatrix}\right] \Diag( \matbf{F}_{k,11}^{-1}\otimes\matbf{I}_{d},\matbf{I}_{md})
\end{equation*} and 
\begin{equation*}
\begin{aligned}
 \vecbf{h}_{k} & = [
(\vecbf{1}_{n+1} \otimes \vecbf{h}_{k}^{\lrangle{0}})^{\top},
- \underline{h}_{k}^{\lrangle{1}} \vecbf{1}_{nd}^{\top},
\overline{h}_{k}^{\lrangle{1}} \vecbf{1}_{nd}^{\top},
\cdots, \\
& \quad \quad - \underline{h}_{k}^{\lrangle{m-1}} \vecbf{1}_{(n-m+2)d}^{\top},
\overline{h}_{k}^{\lrangle{m-1}} \vecbf{1}_{(n-m+2)d}^{\top}	 
]^{\top}.
\end{aligned}
\end{equation*}

\subsection{The overall problem and algorithm}
With the objective function and constraints constructed above, we can formulate the trajectory generation problem as the following multi-stage constrained optimal control problem:
\begin{equation}
\label{eq:prob}
\begin{aligned}
\min_{\myset{U}} &~~ \eqref{eq:obj1} \\
\mathrm{s.t.} &~~\eqref{eq:state}, \eqref{eq:cons}, \vecbf{x}_{0}~\mathrm{is~given}.
\end{aligned}
\end{equation}

It should be noted that the constraint and objective function in \eqref{eq:prob} are respectively in linear and quadratic forms, in fact, all the matrices involved are functions of $t_{k}$ and hence they are not linear-quadratic. Therefore, we shall solve the formulated problem via IPDDP \cite{pavlov2021interior}, which generalizes the traditional unconstrained DDP to a constrained version. Similar to DDP, the IPDDP algorithm iterates between backward pass, which computes control inputs to
minimize a quadratic approximation of cost in the vicinity of the current trajectory, and forward pass, which updates the current trajectory to a new one, until some terminating condition is triggered.

Define $V_{N} = \|\vecbf{x}_{N} - \vecbf{x}_{N,\mathrm{g}}\|_{\matbf{Q}_{N}}^{2}$. Let $\ell_{k} = \|\vecbf{x}_{k} - \vecbf{x}_{k,\mathrm{g}}\|_{\matbf{Q}_{k}}^{2} + \|\vecbf{u}_{k}\|_{\matbf{R}_{k}}^{2} + \gbf{\lambda}_{k}^{\top} \vecbf{g}_{k}$, and $L_{k} = \ell_{k} + V_{k+1}$ for all $k \in \myset{I}_{N}$, where $\gbf{\lambda}_{k} \succeq \vecbf{0}$ is the Lagrange multiplier. In the following, we shall omit subscript $k$ and denote $(\cdot)_{k+1}$ by $(\cdot)^{+}$ for clarity. Letting $(\cdot)_{\vecbf{a}} := \pdv{(\cdot)}{\vecbf{a}}$ and $(\cdot)_{\vecbf{ab}} := \pdv{(\cdot)}{\vecbf{a}}{\vecbf{b}}$, one has the following iteration which resembles traditional DDP:
\begin{equation}
\label{eq:bwd1}
	\begin{aligned}
	 L_{\gbf{\lambda}} &= \vecbf{g}, L_{\gbf{\lambda} \vecbf{x}} = \vecbf{g}_{\vecbf{x}}, L_{\gbf{\lambda}\vecbf{u}} = \vecbf{g}_{\vecbf{u}}, L_{\gbf{\lambda \lambda}} = \vecbf{0},  \\
	 L_{\vecbf{x}} & = \ell_{\vecbf{x}} + \vecbf{f}_{\vecbf{x}}^{\top}V_{\vecbf{x}}^{+}, L_{\vecbf{u}} = \ell_{\vecbf{u}} + \vecbf{f}_{\vecbf{u}}^{\top}V_{\vecbf{u}}^{+}, \\	
	 L_{\vecbf{xx}} & = \ell_{\vecbf{xx}} + \vecbf{f}_{\vecbf{x}}^{\top}V_{\vecbf{xx}}^{+}\vecbf{f}_{\vecbf{x}} + V_{\vecbf{x}}^{+} \odot \vecbf{f}_{\vecbf{xx}}, \\
	 L_{\vecbf{ux}} & = \ell_{\vecbf{ux}} + \vecbf{f}_{\vecbf{u}}^{\top}V_{\vecbf{ux}}^{+}\vecbf{f}_{\vecbf{x}}+ V_{\vecbf{x}}^{+} \odot \vecbf{f}_{\vecbf{ux}}, \\	  
	 L_{\vecbf{uu}} & = \ell_{\vecbf{uu}} + \vecbf{f}_{\vecbf{u}}^{\top}V_{\vecbf{uu}}^{+}\vecbf{f}_{\vecbf{u}}+ V_{\vecbf{x}}^{+} \odot \vecbf{f}_{\vecbf{uu}}, \\	
	\end{aligned}
\end{equation}
where $\odot$ denotes the tensor contraction and $\vecbf{f}$ is defined in \eqref{eq:state}. In addition to the decision variable $\delta \vecbf{u}_{k}$ (i.e., the first order variation of $\vecbf{u}_{k}$), the primal-dual version of DDP introduces the additional decision variable $\delta \gbf{\lambda}_{k}$. Both $\delta \vecbf{u}_{k}$ and $\delta \gbf{\lambda}_{k}$ can be represented in an affine form of $\delta \vecbf{x}_{k}$ by solving the optimality conditions of an approximation of constrained Bellman equation (see $(6)$ in \cite{pavlov2021interior}), i.e.,   
\begin{equation}
\label{eq:control_gain}
	\begin{aligned}
	\begin{bmatrix}
	\vecbf{k}^{\vecbf{u}} & \vecbf{K}^{\vecbf{u}} \\
	\vecbf{k}^{\gbf{\lambda}} & \vecbf{K}^{\gbf{\lambda}}
	\end{bmatrix} = - 	\begin{bmatrix}
	L_{\vecbf{uu}} & L_{\vecbf{u}\gbf{\lambda}} \\
	\gbf{\Lambda}L_{\gbf{\lambda}\vecbf{u}} & \matbf{G}
	\end{bmatrix}^{-1}
		\begin{bmatrix}
	L_{\vecbf{u}} & L_{\vecbf{ux}} \\
	\vecbf{r} & \gbf{\Lambda} L_{\gbf{\lambda}\vecbf{x}}
	\end{bmatrix},
	\end{aligned}
\end{equation}
where  $\vecbf{k}^{\vecbf{u}}$, $\vecbf{K}^{\vecbf{u}}$, $\vecbf{k}^{\gbf{\lambda}}$ and  $\vecbf{K}^{\gbf{\lambda}} $ are the update gains for $\delta \vecbf{u}_{k}$ and $\delta \gbf{\lambda}_{k}$, respectively; $\vecbf{r} = \gbf{\Lambda}\vecbf{g} + \mu\vecbf{1}$, $\matbf{G} = \Diag(\vecbf{g})$, $\gbf{\Lambda} = \Diag(\gbf{\lambda})$ and $\mu$ is a perturbation constant. 

With these at hand, the gradient and Hessian of $V_{k}$ can be updated by 
\begin{equation}
\label{eq:v_upd}
\begin{aligned}
	V_{\vecbf{x}} & = \hat{L}_{\vecbf{x}} + \hat{L}_{\vecbf{xu}} \vecbf{k}^{\vecbf{u}}, \\
	V_{\vecbf{xx}} & = \hat{L}_{\vecbf{xx}} + \hat{L}_{\vecbf{xu}} \vecbf{K}^{\vecbf{u}},
\end{aligned}
\end{equation} with the following equation:
\begin{equation}
\label{eq:aux_mat}
	\begin{aligned}
		\hat{L}_{\vecbf{x}} & = L_{\vecbf{x}} - L_{\vecbf{x}\gbf{\lambda}}\matbf{G}^{-1} \vecbf{r}, \\
		\hat{L}_{\vecbf{u}} & = L_{\vecbf{u}} - L_{\vecbf{u}\gbf{\lambda}}\matbf{G}^{-1} \vecbf{r}, \\
		\hat{L}_{\vecbf{xx}} & = L_{\vecbf{xx}} - L_{\vecbf{x}\gbf{\lambda}}\matbf{G}^{-1}\gbf{\Lambda} L_{\gbf{\lambda}\vecbf{x}}, \\		
		\hat{L}_{\vecbf{ux}} & = L_{\vecbf{ux}} - L_{\vecbf{u}\gbf{\lambda}}\matbf{G}^{-1}\gbf{\Lambda} L_{\gbf{\lambda}\vecbf{x}}, \\	
		\hat{L}_{\vecbf{uu}} & = L_{\vecbf{uu}} - L_{\vecbf{u}\gbf{\lambda}}\matbf{G}^{-1}\gbf{\Lambda} L_{\gbf{\lambda}\vecbf{u}},			
	\end{aligned}
\end{equation}
where  $\delta \gbf{\lambda}_{k}$ was eliminated. 

The forward pass resembles the traditional DDP while updating both the control input $\vecbf{u}_{k}$ and the dual variable $\gbf{\lambda}_{k}$:
\begin{equation}
\label{eq:fwd}
	\begin{aligned}
		\vecbf{x}_{0}^{\dagger} & = \vecbf{x}_{0}, \\
		\vecbf{u}_{k}^{\dagger} & = \vecbf{u}_{k} + \vecbf{k}^{\vecbf{u}}_{k} + \matbf{K}^{\vecbf{u}}_{k} (\vecbf{x}_{k}^{\dagger} - \vecbf{x}_{k}), \\
		\gbf{\lambda}_{k}^{\dagger} & = \gbf{\lambda}_{k} + \vecbf{k}^{\gbf{\lambda}}_{k} + \matbf{K}^{\gbf{\lambda}}_{k} (\vecbf{x}_{k}^{\dagger} - \vecbf{x}_{k}), \\
		\vecbf{x}_{k+1}^{\dagger} & = \vecbf{f}(\vecbf{x}_{k}^{\dagger}, \vecbf{u}_{k}^{\dagger}),		 
	\end{aligned}
\end{equation}
where $\vecbf{x}^{\dagger}$, $\vecbf{u}^{\dagger}$, and $\gbf{\lambda}^{\dagger}$ denote the new state, control variable, and dual variable, respectively.

Algorithm \ref{alg:ipddp} presents the entire algorithm for generating a polynomial trajectory using IPDDP. It should be noted that the presented algorithm requires a feasible initial trajectory. This can be generated from a similar adaptation of the Infeasible-IPDDP algorithm in \cite{pavlov2021interior}. 
 
\begin{algorithm}
	\caption{DIRECT}
	\label{alg:ipddp}
	\begin{algorithmic}
		\REQUIRE $\vecbf{x}_{k, \mathrm{g}}$, $\matbf{Q}_{k}$, $\eta_{m, k}$,  $w_{k}$, $\vecbf{x}_{k, \mathrm{g}}$, $k \in \myset{I}_{N}$, $\matbf{Q}_{N}$, a feasible trajectory with $\myset{C}$ and $\myset{T}$.
		\ENSURE $\myset{C}^{\ast}$, $\myset{T}^{\ast}$.	
		\STATE Convert $\myset{C}$, $\myset{T}$ to $\myset{U}$
		\FOR{$\mathrm{iter} = 1:\mathrm{max\_iteration}$}
		\STATE Set $V_{N}$, $V_{\vecbf{x},N}$, $V_{\vecbf{xx},N}$
		\FOR{$k = N-1, \dots, 0$} 
		\STATE Evaluate \eqref{eq:bwd1}
		\STATE Compute control inputs and auxiliary matrices from \eqref{eq:control_gain} and \eqref{eq:aux_mat}
		\STATE Update $V_{\vecbf{x}}$, $V_{\vecbf{xx}}$ according to \eqref{eq:v_upd}
		\ENDFOR
		\STATE Set $\vecbf{x}_{0}^{\dagger} = \vecbf{x}_{0}$
		\FOR{$k = 0, \dots, N-1$} 
		\STATE Update the control variable $\vecbf{u}_{k}^{\dagger}$, multiplier $\gbf{\lambda}_{k}^{\dagger}$ and next state $\vecbf{x}_{k+1}^{\dagger}$ according to \eqref{eq:fwd}
		\ENDFOR		
		\IF{termination condition is satisfied} 
		\STATE break
		\ENDIF
		\ENDFOR
		\STATE Extract $\myset{U}^{\ast}$ and convert back to $\myset{C}^{\ast}$ and $\myset{T}^{\ast}$
	\end{algorithmic}
\end{algorithm}
\section{Specializations}
\label{sec:spec}
In last section, we have introduced an algorithm for jointly optimizing energy and time subject to safety and dynamical feasibility constraints. In that algorithm, part of the polynomial coefficient $\matbf{C}_{k,2}$ (i.e., $\vecbf{v}_{k}$) and allocated time $t_{k}$ for segment $k$ are both regarded as control inputs in the formulation \eqref{eq:prob}. In this section, we shall show that this framework can be reduced to two other types of simpler problems.

\subsection{Fixed time allocation}
\label{subsec:fixedtime}
 It can be easily found that the above algorithm can be used for constrained polynomial trajectory generation with preallocated time. In particular, letting $t_{k} = \const$ and redefining $\vecbf{u}_{k} = \vecbf{v}_{k}$ for all $k \in \myset{I}_{N}$, $\eqref{eq:state}$ turns to be a linear-time-varying system and $\eqref{eq:cons}$ is a linear-time-varying constraint w.r.t. 
$[\vecbf{x}_{k}^{\top}, \vecbf{v}_{k}^{\top}]^{\top}$\footnote{The ``time-varying" is in the sense of $k$.}. The problem can be formally written as:
\begin{equation}
\label{eq:prob_reduced1}
\begin{aligned}
\min_{\myset{V}} &~~ J(\myset{V}) =\sum_{k=0}^{N-1}(  \|\vecbf{x}_{k} - \vecbf{x}_{k,\mathrm{g}}\|_{\matbf{Q}_{k}}^{2} + \|\vecbf{v}_{k}\|_{\matbf{R}_{k}}^{2}) \\
&\quad \quad \quad \quad \quad + \|\vecbf{x}_{N} - \vecbf{x}_{N,\mathrm{g}}\|_{\matbf{Q}_{N}}^{2} \\
\mathrm{s.t.} &~~\vecbf{x}_{k+1} = \matbf{A}_{k} \vecbf{x}_{k} + \matbf{B}_{k}\vecbf{v}_{k}, \\
&~~\matbf{H}_{k}[\vecbf{x}_{k}^{\top}, \vecbf{v}_{k}^{\top}]^{\top} \preceq \vecbf{h}_{k},  \\
&~~\vecbf{x}_{0}~\mathrm{is~given},
\end{aligned}
\end{equation}
where $\myset{V} := \{\vecbf{v}_{k}\}_{k \in \myset{I}_{N}}$. In this case, Algorithm \ref{alg:ipddp} can still be used for solving such type of problem\footnote{It is also a typical linear-time-varying model predictive control (MPC) problem, which can be solved by existing MPC algorithms.}, with all $\ell$ and $\vecbf{f}$ related terms in \eqref{eq:bwd1} being precomputed and saved during the initialization stage instead of being reevaluated at each iteration.

\subsection{Fixed time allocation and no constraint}
If neither safety nor dynamical feasibility constraint is considered, one has  the following problem:

\begin{equation}
\label{eq:prob_reduced2}
\begin{aligned}
\min_{\myset{V}} &~~ J(\myset{V}) =\sum_{k=0}^{N-1}(  \|\vecbf{x}_{k} - \vecbf{x}_{k,\mathrm{g}}\|_{\matbf{Q}_{k}}^{2} + \|\vecbf{v}_{k}\|_{\matbf{R}_{k}}^{2}) \\
&\quad \quad \quad \quad \quad + \|\vecbf{x}_{N} - \vecbf{x}_{N,\mathrm{g}}\|_{\matbf{Q}_{N}}^{2} \\
\mathrm{s.t.} &~~\vecbf{x}_{k+1} = \matbf{A}_{k} \vecbf{x}_{k} + \matbf{B}_{k}\vecbf{v}_{k}, \\
&~~\vecbf{x}_{0}~\mathrm{is~given}.
\end{aligned}
\end{equation}
We show that \eqref{eq:prob_reduced2} is indeed a LQT problem by proving that $\matbf{R}_{k} \succ \matbf{0}$. Firstly, the positive semi-definiteness $\matbf{R}_{k} = \eta_{m, k} \int_{0}^{t_{k}} \matbf{R}_{\vecbf{v}}(t) \dd t$ can be easily established by the linearity of integral. Secondly, denoting the last $m$ elements of $\vecbf{b}^{(m)}(t)$ by $\gbf{\sigma}(t)$, one has 
\begin{equation*}
	\begin{aligned}
	 \gbf{\sigma}(t) & = [\frac{m!}{0!}t^{0}, \cdots, \frac{n!}{(m-1)!}t^{m-1}]^{\top}  \\
	  & = \Diag( \{(m+i)!\}_{i=0}^{m-1}) [\frac{1}{0!}t^{0}, \cdots, \frac{1}{(m-1)!}t^{m-1}]^{\top} \\
	  & =: \gbf{\Psi} \gbf{\tau}. 
	\end{aligned}
\end{equation*}
On the other hand, $\gbf{\tau}$ can be written as $\gbf{\tau} = e^{\underline{\matbf{A}}^{\top}t}\vecbf{e}_{1}$ with $\underline{\matbf{A}} = \left[\begin{smallmatrix}
\vecbf{0}_{m-1} & \matbf{I}_{m-1} \\
 0 & \vecbf{0}_{m-1}^{\top}
\end{smallmatrix}\right]$ and $\vecbf{e}_{1} = [1, \vecbf{0}_{m-1}^{\top}]^{\top}$. Since the pair $(\underline{\matbf{A}}, \vecbf{e}_{1}^{\top})$ is observable, the matrix $\int_{0}^{\tau}e^{\underline{\matbf{A}}^{\top}t}\vecbf{e}_{1} \vecbf{e}_{1}^{\top}e^{\underline{\matbf{A}}t} \dd t$ is nonsingular for any $\tau > 0$. Therefore,
\begin{equation*}
 \int_{0}^{t_{k}} \matbf{R}_{\vecbf{v}}(t) \dd t =  \int_{0}^{t_{k}} \gbf{\Psi} \gbf{\tau} \gbf{\tau}^{\top} \gbf{\Psi}^{\top} \dd t =  \gbf{\Psi} \left(\int_{0}^{t_{k}} \gbf{\tau} \gbf{\tau}^{\top}  \dd t \right)\gbf{\Psi}^{\top}
\end{equation*} is positive definite and this shows that \eqref{eq:prob_reduced2} is an LQT problem.

Let $\matbf{P}_{N} = \matbf{Q}_{N}$. The backward iteration for the resulting LQT problem reduces to
\begin{equation}
\label{eq:lqt_bwd}
\begin{aligned}
\matbf{K}_{k} & = (\matbf{B}_{k}^{\top}\matbf{P}_{k+1}\matbf{B}_{k} + \matbf{R}_{k})^{-1}\matbf{B}_{k}\matbf{P}_{k+1}\matbf{A}_{k}, \\
	\matbf{P}_{k} & = \matbf{A}_{k}^{\top}\matbf{P}_{k+1}(\matbf{A}_{k} - \matbf{B}_{k}\matbf{K}_{k}) + \matbf{Q}_{k}, \\
	\vecbf{q}_{k} & = (\matbf{A}_{k} - \matbf{B}_{k}\matbf{K}_{k})^{\top}\vecbf{q}_{k+1} + \matbf{Q}_{k}\vecbf{x}_{k,\mathrm{g}}, \\
	\vecbf{k}_{k} & = (\matbf{B}_{k}^{\top}\matbf{P}_{k+1}\matbf{B}_{k} + \matbf{R}_{k})^{-1}\matbf{B}_{k}\vecbf{q}_{k+1},
\end{aligned}
\end{equation}
where $\matbf{K}_{k}$ and $\vecbf{k}_{k}$ are the parameters for expressing the new control inputs as a linear form of current state (resembles $\vecbf{k}^{\vecbf{u}}$ and $\vecbf{K}^{\vecbf{u}}$ in \eqref{eq:control_gain}). The forward pass reads:
\begin{equation}
\label{eq:lqt_fwd}
	\begin{aligned}
	 \vecbf{u}_{k} & = - \matbf{K}_{k}\vecbf{x}_{k} + \vecbf{k}_{k}, \\
	 \vecbf{x}_{k+1} & = \matbf{A}_{k}\vecbf{x}_{k} + \matbf{B}_{k}\vecbf{u}_{k},
	\end{aligned}
\end{equation}
which iteratively computes the control input and next state.

Algorithm \ref{alg:lqt} summarizes the above iteration and this gives the analytical solution to the unconstrained fixed-time polynomial trajectory generation problem \eqref{eq:prob_reduced2}. It should be noted that Algorithm \ref{alg:lqt} can be regarded as a simplified version of Algorithm \ref{alg:ipddp} with exact one iteration.

\begin{algorithm}
	\caption{LQT-based trajectory generation.}
	\label{alg:lqt}
	\begin{algorithmic}
		\REQUIRE $\vecbf{x}_{k, \mathrm{g}}$, $\matbf{Q}_{k}$, $\eta_{m, k}$ ($k \in \myset{I}_{N}$), $\matbf{Q}_{N}$, $\myset{T}$.
		\ENSURE $\myset{C}^{\ast}$.
		\STATE Set $\matbf{P}_{N} = \matbf{Q}_{N}$, $\vecbf{q}_{N} = \vecbf{x}_{N,\mathrm{g}}$
		\FOR{$k = N-1, \dots, 0$} 
		\STATE Compute $\matbf{K}_{k}$, $\vecbf{k}_{k}$ according to \eqref{eq:lqt_bwd}
		\ENDFOR
		\FOR{$k = 0, \dots, N-1$} 
		\STATE Compute the control variable $\vecbf{u}_{k}$ and next state $\vecbf{x}_{k+1}$ according to \eqref{eq:lqt_fwd}
		\ENDFOR		
		\STATE Extract $\myset{U}$ and convert it to $\myset{C}^{\ast}$
	\end{algorithmic}
\end{algorithm}

\section{Simulation and Experimental Results}
\label{sec:results}
In this section, we shall demonstrate the efficacy of our proposed framework by comparing with state-of-the-art methods as well as performing real-world experiments. 
\subsection{Benchmarks for trajectory generation methods}
\begin{figure}
	\centering
	\includegraphics[scale = 0.63,trim={0.5cm 0cm 0cm 0cm},clip]{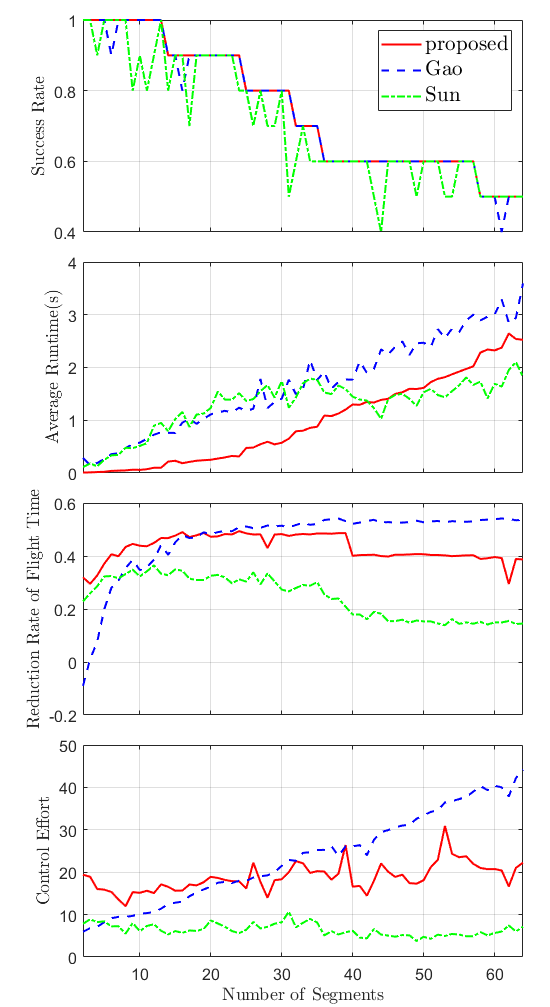}
	\caption{Success rate, average runtime, flight time reduction rate and control effort of the three algorithms. ``Gao'' and ``Sun'' refer to the algorithms proposed in \cite{gao2020teach} and \cite{sun2021fast}, respectively.}
	\label{fig:combined_simulation}
\end{figure}

For trajectory generation with time optimization as well as safety and dynamical feasibility constraints, we compare our algorithm with the alternating spatial-temporal optimization\cite{gao2020teach} and bi-level optimization with analytic gradient\cite{sun2021fast}.\footnote{According to Fig. 9 of \cite{WANG2021GCOPTER}, these two algorithms are the state-of-the-art methods for generating a safe and dynamical feasible trajectory with time optimization. The code for \cite{WANG2021GCOPTER} is currently unavailable.} In \cite{gao2020teach} the energy and time costs are separately optimized while in \cite{WANG2021GCOPTER} the overall cost function takes a similar form as ours, but without the terminal cost and the time cost takes the form of $w_{k}t_k$. The benchmark is conducted using an i7-8550U CPU.
Our algorithm is implemented in three steps: $1)$ initialize our algorithm with Infeasible-IPDDP \cite{pavlov2021interior} and zero initial condition to find a feasible solution with $w_{k} = 1$ and $\matbf{Q}_{N} = \matbf{I}$; $2)$ feed the obtained solution into Algorithm \ref{alg:ipddp} with $w_{k} = 20$ and $\matbf{Q}_{N} = 100\matbf{I}$ to jointly optimize energy and time; and $3)$ feed the obtained solution into Algorithm \ref{alg:ipddp} with fixed time allocation setting (see Section \ref{subsec:fixedtime}). The algorithm in \cite{gao2020teach} is implemented as it is (the coefficient of time is set as $10$). The algorithm in \cite{sun2021fast} is implemented in C++ with MOSEK\footnote{\url{https://www.mosek.com/}} solver for the inner loop QP and backtracking line search for the outer loop and the coefficient of time is set as $20$. We randomly generate $10$ paths with $64$ safety flight corridors (the facet number ranges from $6$ to $117$, and is $28$ on average) and extract $2 \sim 64$ of them to implement these three algorithms. The initial allocation time is set by the subroutine in \cite{gao2020teach}. We consider the problem of generating minimum jerk trajectory with $n=5$, where the maximal velocity and acceleration are set as $2 \mathrm{m/s}$ and $2 \mathrm{m/s^{2}}$, respectively.  The success rate versus number of segments is shown in Fig. \ref{fig:combined_simulation}. The success rate is dropping in general with the increasing number of segments. The success rate of our algorithm is marginally higher than Gao's in some cases since the terminal constraint is softened in our formulation. As a result of terminal constraints and enforcing the dynamical feasibility via control points in inner QP, Sun's method has the lowest success rate. 
The overall computation time and the reduction rate of the flight time (computed by $[\sum_{k}^{N-1} (t_{k}^{\mathrm{ini}} - t_{k}^{\ast})]/\sum_{k}^{N-1} t^{\mathrm{ini}}_{k}$, where $t_{k}^{\mathrm{ini}}$ is the initial allocated time for segment $k$ and $t_{k}^{\ast}$ is the optimized allocated time computed by each algorithm) and the control effort are also compared in Fig. \ref{fig:combined_simulation}. It can be found that for middle number of segments ($N\leq40$), our algorithm consumes less time than the other two algorithms. For $N>40$, our algorithm still consumes less time than Gao's while more than Sun's. The reason is that the outer loop of Sun's method terminates in only a few iterations and gets stuck in local minimum in these cases, which can be seen from the relatively low reduction rate in Fig. \ref{fig:combined_simulation}. 
In addition, there is a trade-off between flight time and control effort: a shorter flight time (i.e., higher reduction rate) corresponds to a higher control effort. 
Unlike Gao's method which has increasing reduction rate and control effort w.r.t. $N$, our algorithm maintains relatively consistent reduction rate and level of control energy. 
Since our formulation also incorporates the terminal cost term (see \eqref{eq:prob}), the change of this term affects the other two terms and hence causes the fluctuation of control energy in our approach, as observed in Fig. \ref{fig:combined_simulation}.
Overall, for the application of trajectory generation with $N \in [10, 30]$, our algorithm achieves comparable results with Gao's method while consumes less runtime.

\begin{figure}
	\centering
	\includegraphics[scale = 0.6,trim={0.1cm 0cm 1cm 0cm},clip]{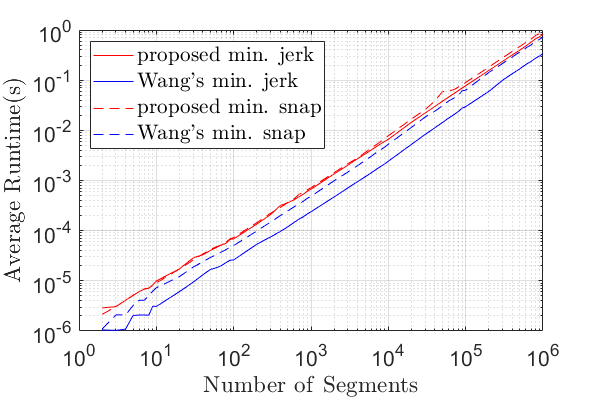}
	\caption{Comparison of large scale trajectory generation.}
	\label{fig:large_scale}
\end{figure}

For trajectory generation with fixed time allocation and without constraint, we compare our Algorithm \ref{alg:lqt} with \cite{wang2020generating}. Our parameters are set as  $\eta_{m, k} = 10^{-5}$, $\matbf{Q}_{k} = 100\matbf{I}$, $k \in \myset{I}_{N}$, $\matbf{Q}_{N} = 100\matbf{I}$. We randomly sample $100$ times for each number of waypoints and run $4$ algorithms ($2$ for minimum jerk trajectory and $2$ for minimum snap trajectory). The average runtime versus number of segments is shown in Fig. \ref{fig:large_scale}. 
It can be found that all the algorithms enjoy linear complexity w.r.t. number of segments. Although our algorithm for minimum jerk trajectory consumes around three times the time of \cite{wang2020generating}, the gap for minimum snap trajectory is much narrower(around 1.5 times).
Furthermore, based on the results reported in \cite{wang2020generating},  our algorithm is faster than the optimized version of algorithm in \cite{burke2020generating}, which is one order of magnitude slower than \cite{wang2020generating} (we do not plot the result of \cite{burke2020generating} as the code is not available.)

\subsection{Flight experiments}
\begin{figure}
	\centering
	\includegraphics[scale = 0.6,trim={0.9cm 0cm 0cm 0cm},clip]{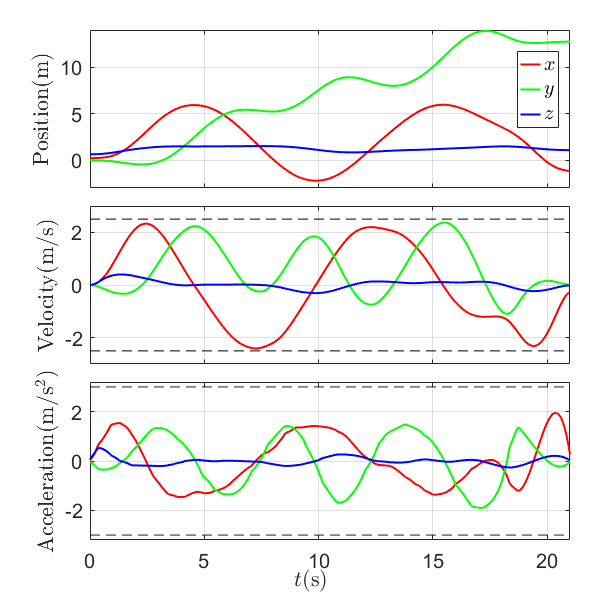}
	\caption{The trajectory command generated from our proposed method.}
	\label{fig:exp_posvelacc}
\end{figure}

The proposed framework is implemented in ROS and we conduct multiple flight experiments of a small quadrotor executing the trajectories generated using our method.
We use stereo camera fusion with IMU\cite{qin2019a} for localization and state estimation of the quadrotor in the environment and use the robust and perfect tracking (RPT) controller \cite{cai2011unmanned} for trajectory tracking.
Firstly, we conduct a set of experiments where a long trajectory is generated in a pre-built map consisting of two rooms and multiple obstacles, as shown in Fig. \ref{fig:longtraj}.
The trajectory is generated by the i7 onboard computer and executed by the quadrotor immediately.
We use multiple settings of dynamical feasibility. The trajectory with maximal velocity and acceleration set as $2.5 \mathrm{m/s}$ and $3\mathrm{m/s^{2}}$ is shown in Fig. \ref{fig:exp_posvelacc}.
In the second set of experiments the onboard computer generates safe and feasible trajectories online to reach target points set by a commander remotely.
The video of the experiments can be found online\footnote{\url{https://youtu.be/BM8_ABM_2VM}} or in the supplemental material.  

\color{black}
\section{Conclusion and Future Work}
\label{sec:conclu}
In this paper, we have proposed to use a state-space equation to characterize a piecewise-polynomial trajectory such that the trajectory generation problem (either with or without time optimization and constraints) can be subsequently solved via DDP algorithm. 
The proposed framework is implemented in C++ without the use of solvers, 
and the efficiency of our proposed framework was demonstrated by comparing with state-of-the-art methods as well as real-world experiments. 

The proposed framework opens up opportunity for further developments.
For example, the joint energy-time optimization can be extended to a multi-robot scenario with non-collision constraints.
The fixed-time optimization can be used for online trajectory replanning due to its fast computation.
\color{black}
\bibliographystyle{ieeetr}        
\bibliography{main}

\begin{thebibliography}{10}

\bibitem{hu2013cooperative}
J.~Hu, J.~Xu, and L.~Xie, ``Cooperative search and exploration in robotic
  networks,'' {\em Unmanned Systems}, vol.~1, no.~1, pp.~121--142, 2013.

\bibitem{ogren2004cooperative}
P.~Ogren, E.~Fiorelli, and N.~E. Leonard, ``Cooperative control of mobile
  sensor networks: Adaptive gradient climbing in a distributed environment,''
  {\em IEEE Transactions on Automatic Control}, vol.~49, no.~8, pp.~1292--1302,
  2004.

\bibitem{yue2019multilevel}
Y.~Yue, C.~Yang, Y.~Wang, P.~C.~N. Senarathne, J.~Zhang, M.~Wen, and D.~Wang,
  ``A multilevel fusion system for multirobot 3-d mapping using heterogeneous
  sensors,'' {\em IEEE Systems Journal}, published online, 2019.

\bibitem{mellinger2011minimum}
D.~Mellinger and V.~Kumar, ``Minimum snap trajectory generation and control for
  quadrotors,'' in {\em 2011 IEEE International Conference on Robotics and
  Automation}, pp.~2520--2525, IEEE, 2011.

\bibitem{bry2015aggressive}
A.~Bry, C.~Richter, A.~Bachrach, and N.~Roy, ``Aggressive flight of fixed-wing
  and quadrotor aircraft in dense indoor environments,'' {\em The International
  Journal of Robotics Research}, vol.~34, no.~7, pp.~969--1002, 2015.

\bibitem{gao2020teach}
F.~Gao, L.~Wang, B.~Zhou, X.~Zhou, J.~Pan, and S.~Shen, ``Teach-repeat-replan:
  A complete and robust system for aggressive flight in complex environments,''
  {\em IEEE Transactions on Robotics}, vol.~36, no.~5, pp.~1526--1545, 2020.

\bibitem{pavlov2021interior}
A.~Pavlov, I.~Shames, and C.~Manzie, ``Interior point differential dynamic
  programming,'' {\em IEEE Transactions on Control Systems Technology}, early
  access, 2021.

\bibitem{liu2017planning}
S.~Liu, M.~Watterson, K.~Mohta, K.~Sun, S.~Bhattacharya, C.~J. Taylor, and
  V.~Kumar, ``Planning dynamically feasible trajectories for quadrotors using
  safe flight corridors in 3-d complex environments,'' {\em IEEE Robotics and
  Automation Letters}, vol.~2, no.~3, pp.~1688--1695, 2017.

\bibitem{sun2021fast}
W.~Sun, G.~Tang, and K.~Hauser, ``Fast uav trajectory optimization using
  bilevel optimization with analytical gradients,'' {\em IEEE Transactions on
  Robotics}, early access, 2021.

\bibitem{wang2020alternating}
Z.~Wang, X.~Zhou, C.~Xu, J.~Chu, and F.~Gao, ``Alternating minimization based
  trajectory generation for quadrotor aggressive flight,'' {\em IEEE Robotics
  and Automation Letters}, vol.~5, no.~3, pp.~4836--4843, 2020.

\bibitem{WANG2021GCOPTER}
Z.~Wang, X.~Zhou, C.~Xu, and F.~Gao, ``Geometrically constrained trajectory
  optimization for multicopters,'' {\em arXiv preprint arXiv:2103.00190}, 2021.

\bibitem{burke2020generating}
D.~Burke, A.~Chapman, and I.~Shames, ``Generating minimum-snap quadrotor
  trajectories really fast,'' in {\em 2020 IEEE/RSJ International Conference on
  Intelligent Robots and Systems (IROS)}, pp.~1487--1492, IEEE, 2020.

\bibitem{wang2020generating}
Z.~Wang, H.~Ye, C.~Xu, and F.~Gao, ``Generating large-scale trajectories
  efficiently using double descriptions of polynomials,'' {\em arXiv preprint
  arXiv:2011.02662}, 2020.

\bibitem{burke2021fast}
D.~Burke, A.~Chapman, and I.~Shames, ``Fast spline trajectory planning: Minimum
  snap and beyond,'' {\em arXiv preprint arXiv:2105.01788}, 2021.

\bibitem{mayne1966second}
D.~Mayne, ``A second-order gradient method for determining optimal trajectories
  of non-linear discrete-time systems,'' {\em International Journal of
  Control}, vol.~3, no.~1, pp.~85--95, 1966.

\bibitem{de1988differential}
J.~De~O.~Pantoja, ``Differential dynamic programming and newton's method,''
  {\em International Journal of Control}, vol.~47, no.~5, pp.~1539--1553, 1988.

\bibitem{morimoto2003minimax}
J.~Morimoto, G.~Zeglin, and C.~G. Atkeson, ``Minimax differential dynamic
  programming: Application to a biped walking robot,'' in {\em Proceedings 2003
  IEEE/RSJ International Conference on Intelligent Robots and Systems (IROS
  2003)(Cat. No. 03CH37453)}, vol.~2, pp.~1927--1932, IEEE, 2003.

\bibitem{kumar2016optimal}
V.~Kumar, E.~Todorov, and S.~Levine, ``Optimal control with learned local
  models: Application to dexterous manipulation,'' in {\em 2016 IEEE
  International Conference on Robotics and Automation (ICRA)}, pp.~378--383,
  IEEE, 2016.

\bibitem{plancher2017constrained}
B.~Plancher, Z.~Manchester, and S.~Kuindersma, ``Constrained unscented dynamic
  programming,'' in {\em 2017 IEEE/RSJ International Conference on Intelligent
  Robots and Systems (IROS)}, pp.~5674--5680, IEEE, 2017.

\bibitem{xie2017differential}
Z.~Xie, C.~K. Liu, and K.~Hauser, ``Differential dynamic programming with
  nonlinear constraints,'' in {\em 2017 IEEE International Conference on
  Robotics and Automation (ICRA)}, pp.~695--702, IEEE, 2017.

\bibitem{aoyama2020constrained}
Y.~Aoyama, G.~Boutselis, A.~Patel, and E.~A. Theodorou, ``Constrained
  differential dynamic programming revisited,'' {\em arXiv preprint
  arXiv:2005.00985}, 2020.

\bibitem{zhang1997splines}
Z.~Zhang, J.~Tomlinson, and C.~Martin, ``Splines and linear control theory,''
  {\em Acta Applicandae Mathematica}, vol.~49, no.~1, pp.~1--34, 1997.

\bibitem{sun2000control}
S.~Sun, M.~B. Egerstedt, and C.~F. Martin, ``Control theoretic smoothing
  splines,'' {\em IEEE Transactions on Automatic Control}, vol.~45, no.~12,
  pp.~2271--2279, 2000.

\bibitem{tordesillas2020minvo}
J.~Tordesillas and J.~P. How, ``Minvo basis: Finding simplexes with minimum
  volume enclosing polynomial curves,'' {\em arXiv preprint arXiv:2010.10726},
  2020.

\bibitem{qin2019a}
T.~Qin, J.~Pan, S.~Cao, and S.~Shen, ``A general optimization-based framework
  for local odometry estimation with multiple sensors,'' 2019.

\bibitem{cai2011unmanned}
G.~Cai, B.~M. Chen, and T.~H. Lee, {\em Unmanned rotorcraft systems}.
\newblock Springer Science \& Business Media, 2011.

\end{thebibliography}

\end{document}